\documentclass[11pt]{amsart}
\usepackage[T1]{fontenc}
\usepackage[utf8]{inputenc}
\usepackage{lmodern}
\usepackage{microtype}
\usepackage{xcolor}
\usepackage{amsmath,amssymb,amsthm,mathtools}
\usepackage{enumitem}
\usepackage[margin=1in]{geometry}
\usepackage{hyperref}
\hypersetup{colorlinks=true,linkcolor=blue,citecolor=blue,urlcolor=blue}
\usepackage{mathrsfs}

\theoremstyle{plain}
\newtheorem{theorem}{Theorem}[section]
\newtheorem{lemma}[theorem]{Lemma}
\newtheorem{proposition}[theorem]{Proposition}

\theoremstyle{remark}
\newtheorem{remark}[theorem]{Remark}

\DeclareMathOperator{\dist}{dist}
\newcommand{\N}{\mathbb{N}}
\newcommand{\C}{\mathbb{C}}
\numberwithin{equation}{section}

\hypersetup{pdftitle={The Kaplansky--Aupetit conjecture for C*-algebras},pdfauthor={Ilja Gogic,Bojan Kuzma},pdfsubject={Invertibility preservers and Jordan homomorphisms}}

\title[The Kaplansky--Aupetit conjecture for $C^*$-algebras]{The Kaplansky--Aupetit conjecture for $C^*$-algebras}

\author{Ilja Gogi\'{c}}
\address{I.~Gogi\'c, University of Zagreb Faculty of Science, Department of Mathematics, Bijeni\v{c}ka 30, 10000 Zagreb, Croatia}
\email{ilja@math.hr, ilja.gogic@math.pmf.unizg.hr}

\author{Bojan Kuzma}
\address{B.~Kuzma, University of Primorska, Glagolja\v{s}ka 8, SI-6000 Koper, Slovenia, and Institute of Mathematics, Physics, and Mechanics, Jadranska 19, SI-1000 Ljubljana, Slovenia}
\email{bojan.kuzma@upr.si}

\date{\today}



\begin{document}

\begin{abstract}
We settle the Kaplansky--Aupetit conjecture when either algebra is a $C^*$-algebra. More precisely, if $A$ is a unital Banach algebra and $B$ is a unital $C^*$-algebra, or if $A$ is a unital $C^*$-algebra and $B$ is a unital semisimple Banach algebra, then every unital invertibility-preserving linear surjection $\phi: A\to B$ is a Jordan homomorphism.  In fact, it suffices to assume that $\phi$ preserves the invertibility of products of two exponentials, or of individual exponentials when both algebras are $C^*$-algebras.

Claude~Opus~5.5 found the original proof when both $A$ and $B$ are $C^*$-algebras by adapting a resolvent estimate from a recent preprint of Colbrook and Stepaniants, while OpenAI's GPT-6.1 further generalized this to the case where only one algebra is a $C^*$-algebra.
\end{abstract}

\maketitle

\section{Introduction}

The celebrated Gleason--Kahane--\.{Z}elazko theorem~\cite{Gleason1967,KahaneZelazko1968,Zelazko1968} states that a linear functional on a unital complex Banach algebra is multiplicative if it takes the identity to $1$ and does not vanish
on invertible elements. This led Kaplansky~\cite[Section~9]{Kaplansky1970} to ask in 1970 whether a similar principle holds for invertibility-preserving linear maps between more general algebras.

To formulate the question more precisely, let $A$ be a nonzero unital Banach algebra with identity~$1_A$ and group of invertible elements $A^\times$. For $x\in A$, we let
\[
    \sigma_A(x):=\{\lambda\in\C:\lambda1_A-x\notin A^\times\} \qquad \text{and}
    \qquad r_A(x):=\max\{|\lambda|:\lambda\in\sigma_A(x)\}
\]
be its spectrum and spectral radius. If $B$ is another unital Banach algebra, then a linear map $\phi\colon A\to B$ is \emph{unital} if $\phi(1_A)=1_B$ and it \emph{preserves invertibility} if $\phi(A^\times)\subseteq B^\times$. For unital maps, this can be  equivalently expressed with the spectrum as $\sigma_B(\phi(x))\subseteq\sigma_A(x)$ for all $x\in A$. The strong form of the \emph{Kaplansky--Aupetit conjecture} asks whether every unital invertibility-preserving linear surjection $\phi\colon A\to B$ between unital semisimple Banach algebras (those with zero Jacobson radical) is a \emph{Jordan homomorphism}, that is, whether it preserves squares ($\phi(x^2)=\phi(x)^2$ for all $x\in A$); see, e.g.,~\cite[Conjecture~1.1]{BresarSemrl2008}. Notice that some form of semisimplicity is essential here, as the unital spectrum-preserving but non-Jordan bijection of the block upper triangular subalgebra of $M_4(\C)$ given by
\[
    \phi\colon\begin{pmatrix}a&b\\0&c\end{pmatrix}\mapsto
    \begin{pmatrix}a&b\\0&c^t\end{pmatrix},
    \qquad a,b,c\in M_2(\C)
\]
clearly demonstrates (see also~\cite[p.~28]{Aupetit1979a}). The converse statement, however, is trivial and holds even without semisimplicity: Every unital Jordan surjection preserves Jordan triple product $xyx$ and as such preserves invertibility because $x$ is invertible if and only if $xyx$ attains $1$ (cf.~also~\cite[Proposition~1.3]{Sourour1996}).

With a harmless normalization  $x\mapsto \phi(1)^{-1}\phi(x)$ we can make any invertibility preserving surjection $\phi\colon A\to B$ into a unital invertibility preserving surjection. If $B$ is semisimple, then every such linear surjection is automatically bounded by~\cite[Theorem~5.5.2]{Aupetit1991}. For maps between von Neumann algebras, Aupetit~\cite[Theorem~1.3,  cf. also Remark~2.7]{Aupetit2000} proved that unital invertibility preserving bijections are Jordan; he obtained the same result even without assuming injectivity if invertibility is preserved in both directions (i.e., $x\in A^\times\iff\phi(x)\in B^\times$ for all $x\in A$). This was extended by  Cui and Hou who, in~\cite[Theorem~2.3]{CuiHou2004}, proved  the same result for unital linear invertibility preserving surjections from a von Neumann algebra onto a unital semisimple Banach algebra. A key step  in all these  proofs is preservation of idempotents, a strategy first recognized by Bre\v{s}ar and \v{S}emrl~\cite{BresarSemrl1998} and subsequently used in related work on linear preservers (see, e.g.,~\cite{AlaminosEtAl2012,Kovacs2005,MathieuSchulz2025}). Further early cases of the conjecture appear in~\cite{Aupetit1979,ChoiEtAl1984,Sourour1996}.

A major breakthrough came in the recent preprint of Colbrook and Stepaniants~\cite[Theorem~1.1]{ColbrookStepaniants2026}, who settled Mathieu and Schick's spectral isometry conjecture~\cite[Conjecture~4.6]{MathieuSchick2002}. They proved that unital linear surjections between unital $C^*$-algebras that preserve spectral radii are Jordan isomorphisms.  In this paper we prove the strong form of the Kaplansky--Aupetit conjecture when at least one of $A$ and $B$ is a $C^*$-algebra. In either case, the unital linear surjection need only send products of two exponentials to invertible elements. Our proof begins with the case of two $C^*$-algebras, where we adapt the resolvent argument of Colbrook and Stepaniants to show that individual exponentials suffice. To pass from exponentials to spectra, recall that the \emph{polynomially convex hull} of a nonempty compact set $K\subseteq\C$ is
\begin{equation}\label{eq:polyhull}
    \widehat K:=
    \left\{\lambda\in\C:
    |p(\lambda)|\le\max_{\zeta\in K}|p(\zeta)|
    \text{ for every }p\in\C[z]\right\},
\end{equation}
where $\C[z]$ is the algebra of complex polynomials in one variable. For $x\in A$, we define its \emph{spectral hull} by
$\widehat{\sigma}_A(x):=\widehat{\sigma_A(x)}$.  If $\lambda\in\C\setminus\widehat{\sigma}_A(x)$, then $\lambda$ lies in the unbounded component of $\C\setminus\sigma_A(x)$, so there is a simple polygonal path in this component starting at $\lambda$ and tending to infinity. Its complement in $\C$ is a simply connected neighborhood of $\sigma_A(x)$ on which $z\mapsto\lambda-z$ has a holomorphic logarithm $f$. The holomorphic functional calculus then gives $e^{f(x)}=\lambda1_A-x$ (see~\cite[proof of Theorem~3.3.6]{Aupetit1991}). Consequently, if a unital linear map $\phi:A\to B$ sends exponentials to invertible elements, the spectral hull of $\phi(x)$ is contained in that of $x$ for every $x\in A$ (Lemma~\ref{lem:exterior}). Our first major result shows that this inclusion already forces a unital linear surjection between $C^*$-algebras to be a Jordan homomorphism.

\begin{theorem}\label{thm:main}
Let $A$ be a unital Banach algebra, let $B$ be a unital semisimple Banach algebra, and let $\phi:A\to B$ be a unital linear surjection.
\begin{enumerate}[label=\textup{(\alph*)},ref=\textup{(\alph*)}]
\item\label{item:main-two} If at least one of $A$ and $B$ is a $C^*$-algebra, then $\phi$ is a Jordan homomorphism if and only if $\phi(e^xe^y)\in B^\times$ for all $x,y\in A$.
\item\label{item:main-single} If both $A$ and $B$ are $C^*$-algebras, then the following conditions are equivalent.
\begin{enumerate}[label=\textup{(\roman*)}]
\item $\phi$ is a Jordan homomorphism.
\item $\phi(e^x)\in B^\times$ for every $x\in A$.
\item The spectral hulls satisfy $\widehat{\sigma}_B(\phi(x))\subseteq\widehat{\sigma}_A(x)$ for all $x\in A$.
\end{enumerate}
\end{enumerate}
\end{theorem}

Aupetit and Mouton~\cite[p.~92]{AupetitMouton1994} asked whether a linear surjection between unital semisimple Banach algebras that preserves spectral hulls must be a Jordan homomorphism.  They showed that every such surjection is unital~\cite[Theorem~3.1(2)]{AupetitMouton1994} and that it is a Jordan homomorphism when $B$ has an essential socle~\cite[Corollary~3.3]{AupetitMouton1994}. Since the spectral hull determines the spectral radius, the spectral isometry theorem of Colbrook and Stepaniants gives an affirmative answer when both algebras are $C^*$-algebras. Under the weaker hypothesis of spectral hull inclusion, Cui and Hou~\cite[Theorem~2.2]{CuiHou2004Spectral} showed that a linear surjection between unital semisimple Banach algebras is unital and preserves idempotents. When the domain is a $C^*$-algebra of real rank zero, the map is also a Jordan homomorphism~\cite[Theorem~3.1]{CuiHou2004Spectral}.

The proof of Theorem~\ref{thm:main}, given in Section~\ref{sec:main}, has two stages. For part~\ref{item:main-single}, hull inclusion makes $\phi$ bounded (Lemma~\ref{lem:exterior}\ref{item:hull}). We pass to the second adjoint $\phi^{**}:A^{**}\to B^{**}$ because on the von Neumann algebra $A^{**}$ it suffices to show that projections have idempotent images (Lemma~\ref{lem:jordan}\ref{item:projection}). For a projection $p\in A^{**}$, we approximate $1-2p$ in the strong$^*$ topology by two nets of unitaries in $A$ with spectra in opposite semicircles. From each point outside $\{-1,1\}$, we can choose a ray that stays away from one of the semicircles. This lets us apply the resolvent estimate adapted from Colbrook and Stepaniants~\cite[proof of Lemma~3.2]{ColbrookStepaniants2026} to the images of the corresponding net (Lemma~\ref{lem:res}). Passing to the limit bounds the resolvent norm of $\phi^{**}(1-2p)$ by a constant times the inverse distance to $\{-1,1\}$. Holomorphic functional calculus then gives $\phi^{**}(1-2p)^2=1$, so $\phi^{**}(p)$ is idempotent (Proposition~\ref{prop:projections}).

For part~\ref{item:main-two}, we multiply the argument of $\phi$ by an exponential on either side and make each resulting map unital. The two-exponential hypothesis ensures that both maps send each $x\in A$ to an element with the same spectral hull as $\phi(x)$. This shows that $\ker\phi$ is a closed two-sided ideal (Lemma~\ref{lem:translates}). Put $C:=A/\ker\phi$ and let $S:C\to B$ be the induced bijection. Products of two exponentials in $C$ lift to such products in $A$, so $S$ satisfies the two-exponential hypothesis. Thus for each exponential $g\in C$, the map $U_g:B\to B$ given by $U_g(S(x)):=S(g)^{-1}S(gx)$ for $x\in C$ is a unital linear bijection preserving spectral hulls (Lemma~\ref{lem:translates}). We apply part~\ref{item:main-single} to $U_g$ when $B$ is a $C^*$-algebra. When $A$ is a $C^*$-algebra, hull preservation by $S$ (Lemma~\ref{lem:periphery}) lets us apply that part to $S^{-1}U_gS$ on the $C^*$-algebra $C$. Differentiating the resulting square identities along exponential paths through $1_C$ gives bounded Jordan derivations, which are derivations by Sinclair's theorem~\cite[Theorem~3.3]{Sinclair1970}. The derivation criterion applies to $S$ in the first case and to $S^{-1}$ in the second (Lemma~\ref{lem:tangent}). Thus $S$, and hence $\phi$, preserves squares in either case.

This paper is organized as follows. Section~\ref{sec:prelim} collects the preliminaries, and Section~\ref{sec:mixed} prepares the reduction to two $C^*$-algebras. The proof of Theorem~\ref{thm:main} is given in Section~\ref{sec:main}.

\subsection*{Note added}
 After the first version of this paper was finished, we learned that Aznay and Chabbabi had posted a preprint on arXiv in which, with the help of OpenAI's ChatGPT, they announce a proof of  the strong form of the Kaplansky--Aupetit conjecture for unital Banach algebras with semisimple codomain~\cite[Proposition~5.2]{AznayChabbabi2026}. With the spectral-radius identities in Lemma~\ref{lem:translates}, their argument also gives the two-exponential criterion in this generality. Both proofs use these identities to show that the kernel is a two-sided ideal. Our proof of that criterion is more restrictive in that we assume that at least one of the algebras $A$ and $B$ is a $C^*$-algebra. However, our techniques are different, and for maps between two $C^*$-algebras we use a resolvent estimate and passage to the bidual to show that it suffices to send individual exponentials to invertible elements or, equivalently, to assume spectral hull inclusion (see Theorem~\ref{thm:main}\ref{item:main-single}).

\section{Preliminaries}\label{sec:prelim}

We collect the facts about spectral hulls, strong$^*$ convergence and Jordan homomorphisms used in the proofs. Throughout, all algebras and linear maps are over $\C$, and unital algebras are assumed to be nonzero. We write $\N$ for the positive integers, $\Lambda$ for the directed index set of a net and $1$ for the identity when the algebra is clear from the context.

An \emph{idempotent} in an algebra $A$ is an element $e\in A$ satisfying $e^2=e$. A self-adjoint idempotent in a $C^*$-algebra is called a \emph{projection}.

Recall that the polynomially convex hull $\widehat K$ of a nonempty compact set $K\subseteq\C$, defined in~\eqref{eq:polyhull}, is obtained by adjoining to $K$ all bounded components of $\C\setminus K$~\cite[Proposition~VII.5.3]{Conway1990}. Thus $\C\setminus\widehat K$ is the unbounded component of $\C\setminus K$.

\begin{lemma}\label{lem:exterior}
Let $A$ and $B$ be unital Banach algebras.
\begin{enumerate}[label=\textup{(\alph*)},ref=\textup{(\alph*)},
    beginpenalty=10000,midpenalty=10000]
\item\label{item:exterior} If $x\in A$ and $\lambda\in\C\setminus\widehat{\sigma}_A(x)$, then $\lambda1_A-x$ is an exponential in $A$.
\item\label{item:hull} If a unital linear map $\phi:A\to B$ sends exponentials to invertible elements, then $\widehat{\sigma}_B(\phi(x))\subseteq\widehat{\sigma}_A(x)$ for all $x\in A$. For a unital linear surjection $\phi:A\to B$ with $B$ semisimple, the hull inclusion alone implies that $\phi$ is bounded.
\end{enumerate}
\end{lemma}

\begin{proof}
Let $x\in A$ and $\lambda\in\C\setminus\widehat{\sigma}_A(x)$, and put $b:=\lambda1_A-x\in A$. Since $\sigma_A(b)=\{\lambda-\zeta:\zeta\in\sigma_A(x)\}$, the map $z\mapsto\lambda-z$, $z\in\C$, carries the unbounded component of $\C\setminus\sigma_A(x)$ onto that of $\C\setminus\sigma_A(b)$. In particular, $0$ belongs to the latter component. Thus $b$ has a logarithm in $A$ by~\cite[Theorem~3.3.6]{Aupetit1991}.

Suppose that $\phi$ sends exponentials to invertible elements. By~\ref{item:exterior}, for every $x\in A$ and
$\lambda\in\C\setminus\widehat{\sigma}_A(x)$ we have $\phi(\lambda1_A-x)=\lambda1_B-\phi(x)\in B^\times$. Thus $\sigma_B(\phi(x))\subseteq\widehat{\sigma}_A(x)$, and taking polynomially convex hulls gives the required inclusion.

Now suppose that $\phi$ is surjective, $B$ is semisimple, and the hull inclusion holds for every $x\in A$. Taking the hull does not change the maximum modulus, so $r_B(\phi(x))\le r_A(x)$ for every $x\in A$. Aupetit's automatic continuity theorem~\cite[Theorem~5.5.2]{Aupetit1991} then gives boundedness.
\end{proof}

The bidual $A^{**}$ of a $C^*$-algebra $A$ is a von Neumann algebra. We identify it with the strong operator closure of the image of $A$ in its universal representation~\cite[Section~10.1]{KadisonRingrose1997b}. The canonical embedding $A\to A^{**}$ is an isometric $*$-homomorphism and is unital when $A$ is unital. By Kaplansky's density theorem~\cite[Theorem~5.3.5]{KadisonRingrose1997a}, every self-adjoint contraction in $A^{**}$ is a strong operator limit of self-adjoint contractions in $A$.

Let $M$ be a von Neumann algebra faithfully and normally represented on a Hilbert space $H$. A norm-bounded net $(x_j)_{j\in\Lambda}$ in $M$ converges to $x\in M$ in the strong$^*$ topology if and only if $x_j\xi\to x\xi$ and $x_j^*\xi\to x^*\xi$ for every $\xi\in H$. For self-adjoint elements, strong operator convergence implies strong$^*$ convergence. On norm-bounded sets the strong$^*$ topology agrees with the $\sigma$-strong$^*$ topology~\cite[Remark~I.2]{Akemann1967}, which is generated by the seminorms $x\mapsto\omega(x^*x+xx^*)^{1/2}$, $x\in M$, as $\omega$ ranges over the positive normal linear functionals on $M$. In particular, on bounded sets the strong$^*$ topology is independent of the faithful normal representation. Multiplication is jointly strong$^*$-continuous on such sets.

\begin{lemma}\label{lem:functional}
Let $(x_j)_{j\in\Lambda}$ be a net of self-adjoint elements of a von Neumann algebra $M$ converging strong$^*$ to a self-adjoint element $x\in M$. Suppose that a compact interval $I\subseteq\mathbb R$ contains the spectrum of $x$ and that of $x_j$ for every $j\in\Lambda$. For every continuous function $f:I\to\C$,
\[
    f(x_j)\longrightarrow f(x)\qquad\text{in the strong$^*$ topology}.
\]
\end{lemma}

\begin{proof}
Choose continuous functions $g,h:\mathbb R\to\mathbb R$ of compact support with $f=g+ih$ on $I$. In a faithful normal
representation of $M$, Kaplansky's lemma~\cite[Lemma~5]{Kaplansky1951} gives strong convergence of $g(x_j)$ to $g(x)$ and of $h(x_j)$ to $h(x)$. Both limits are strong$^*$ because the operators are self-adjoint. Taking the corresponding linear combination gives $f(x_j)\to f(x)$ strong$^*$.
\end{proof}

A bounded linear map $\phi:A\to B$ between Banach spaces extends to its second adjoint $\phi^{**}:A^{**}\to B^{**}$, which is bounded and weak$^*$-continuous. If $A$ and $B$ are $C^*$-algebras, Akemann's theorem~\cite[Corollary~III.10]{Akemann1967} shows that $\phi^{**}$ is strong$^*$-continuous on norm-bounded sets (see also~\cite[Lemma~2.1]{ColbrookStepaniants2026}).

We record two standard facts about Jordan homomorphisms. The first follows from Sourour~\cite[Proposition~1.3]{Sourour1996}, and the second is the projection criterion of Mathieu and Schick~\cite[Lemma~2.1]{MathieuSchick2003}.

\begin{lemma}\label{lem:jordan}
\begin{enumerate}[label=\textup{(\alph*)},ref=\textup{(\alph*)}]
\item\label{item:inverse} A unital Jordan homomorphism between unital Banach algebras preserves invertibility.
\item\label{item:projection} A bounded linear map from a von Neumann algebra into a Banach algebra is a Jordan homomorphism if it sends projections to idempotents.
\end{enumerate}
\end{lemma}

\section{Normalized translates and spectral hulls}\label{sec:mixed}

Let $A$ and $B$ be unital Banach algebras with $B$ semisimple, and let $\phi:A\to B$ be a unital linear surjection such that $\phi(e^xe^y)\in B^\times$ for all $x,y\in A$. For the proof of Theorem~\ref{thm:main}\ref{item:main-two}, we show that $\ker\phi$ is a closed two-sided ideal, give a derivation criterion for preservation of squares, and prove that the induced bijection $A/\ker\phi\to B$ preserves spectral hulls when $A$ is a $C^*$-algebra.

\subsection{Normalized translates and the kernel}

For an exponential $g\in A$, we obtain the \emph{normalized translates} of $\phi$ by multiplying its argument by $g$ and the resulting image by $\phi(g)^{-1}$, both on the left or both on the right. Mathieu and Schulz~\cite[Section~4]{MathieuSchulz2023} observed that products of two exponentials suffice for the Jordan conclusion when the codomain is a $C^*$-algebra with a faithful tracial state. Here we use the spectral Liouville theorem to show that the images of each $x\in A$ under $\phi$ and its two translates have the same spectral hull.

\begin{lemma}\label{lem:translates}
The map $\phi$ is bounded and $\ker\phi$ is a closed two-sided ideal of $A$. For $a,x\in A$ and $g:=e^a\in A$,
\begin{equation}\label{eq:normalized}
    \widehat{\sigma}_B\left(\phi(g)^{-1}\phi(gx)\right)
    =\widehat{\sigma}_B(\phi(x))
    =\widehat{\sigma}_B\left(\phi(xg)\phi(g)^{-1}\right).
\end{equation}
\end{lemma}

\begin{proof}
Taking one of the exponentials to be $1_A$, we see that $\phi$ sends exponentials to invertible elements. Thus $\phi$ is bounded by Lemma~\ref{lem:exterior}\ref{item:hull}.

Fix $a\in A$. For $\zeta\in\C$, put $g_\zeta:=e^{\zeta a}\in A$ and define $F_\zeta,G_\zeta:A\to B$ by
\[
    F_\zeta(x):=\phi(g_\zeta)^{-1}\phi(g_\zeta x),
    \qquad G_\zeta(x):=\phi(xg_\zeta)\phi(g_\zeta)^{-1},
    \qquad x\in A.
\]
Both maps are unital and send exponentials to invertible elements by the two-exponential hypothesis. Fix $x\in A$. By Lemma~\ref{lem:exterior}\ref{item:hull}, the spectra of $F_\zeta(x)$ and $G_\zeta(x)$ lie in the fixed compact set $\widehat{\sigma}_A(x)$ for all $\zeta\in\C$. The maps $\zeta\mapsto F_\zeta(x)$ and $\zeta\mapsto G_\zeta(x)$ are entire, since $\phi$ is bounded and inversion is analytic on $B^\times$. By Aupetit's spectral Liouville theorem~\cite[Theorem~3.4.14]{Aupetit1991}, the spectral hulls of $F_\zeta(x)$ and $G_\zeta(x)$ are independent of $\zeta$. Evaluating at $\zeta=0$ and $\zeta=1$ gives~\eqref{eq:normalized}.

Fix $\zeta\in\C$ and let $T:A\to B$ be either $F_\zeta$ or $G_\zeta$. Each translate is the composition of $\phi$ with invertible multiplication operators on the appropriate side, so $T$ is surjective. For $k\in\ker\phi$ and $x\in A$, taking maximum moduli in~\eqref{eq:normalized} gives
\[
    r_B(T(x)+T(k))
    \stackrel{\eqref{eq:normalized}}{=}r_B(\phi(x+k))
    =r_B(\phi(x))
    \stackrel{\eqref{eq:normalized}}{=}r_B(T(x)).
\]
Since $T$ is surjective, $r_B(b+T(k))=0$ for every quasinilpotent $b\in B$. By Zem\'anek's characterization~\cite[Theorem~5.3.1(iii)]{Aupetit1991} and semisimplicity of $B$, we have $T(k)=0$. Applying this to both translates gives $\phi(e^{\zeta a}k)=\phi(ke^{\zeta a})=0$ for every $\zeta\in\C$. Differentiating at $\zeta=0$ gives $\phi(ak)=\phi(ka)=0$. Thus $\ker\phi$ is a two-sided ideal, and it is closed by boundedness of $\phi$.
\end{proof}

\subsection{A derivation criterion}

In the proof of Theorem~\ref{thm:main}\ref{item:main-two}, we apply the following criterion to the induced bijection $S:A/\ker\phi\to B$ when $B$ is a $C^*$-algebra. When $A$ is a $C^*$-algebra, so is $A/\ker\phi$, and we apply the criterion to $S^{-1}$.

\begin{lemma}\label{lem:tangent}
Let $C$ be an associative complex algebra, let $D$ be a $C^*$-algebra, and let $T:C\to D$ be a linear bijection. For $c\in C$ and $b\in D$, let $L_c:C\to C$ and $L_b:D\to D$ denote left multiplication by $c$ and $b$, respectively. If the operator $\delta_c:D\to D$ given by
\[
    \delta_c:=TL_cT^{-1}-L_{T(c)}
\]
is a derivation for every $c\in C$, then $T$ is a Jordan homomorphism.
\end{lemma}

\begin{proof}
Fix $c\in C$ and put $b:=T(c)\in D$ and $\delta:=\delta_c$. Put $d:=\delta(b)=T(c^2)-b^2\in D$. It suffices to show that $d=0$. The Leibniz rule gives $\delta L_b=L_d+L_b\delta$. Using $L_{c^2}=L_c^2$ and $T(c^2)=b^2+d$, we obtain
\[
    \delta_{c^2}=(L_b+\delta)^2-L_{b^2+d}=2L_b\delta+\delta^2.
\]
For $y\in D$, the Leibniz rule gives $b\delta(by)=b\delta(b)y+b^2\delta(y)$. Thus the terms from $2L_b\delta$ cancel in the Leibniz identity for $\delta_{c^2}$ at $(b,y)$. Expanding $\delta^2(by)$ with the same rule gives
\begin{equation}\label{eq:tangent-zero}
    0=\delta^2(by)-\delta^2(b)y-b\delta^2(y)=2d\delta(y).
\end{equation}
Taking $y=d^*b$ and $y=d^*$ in this identity and using the Leibniz rule for $\delta$ yields
\[
    0\stackrel{\eqref{eq:tangent-zero}}{=}d\delta(d^*b)
    =d\delta(d^*)b+dd^*d
    \stackrel{\eqref{eq:tangent-zero}}{=}dd^*d.
\]
Multiplying on the left by $d^*$ gives $(d^*d)^2=0$. By the $C^*$-identity, $\|d\|^4=\|(d^*d)^2\|=0$, so $d=0$.
\end{proof}

\subsection{Spectral hulls for a \texorpdfstring{$C^*$}{C*}-algebra domain}

When $A$ is a $C^*$-algebra, the following lemma shows that the induced bijection $A/\ker\phi\to B$ preserves spectral hulls. The key step is to show that the boundary of each element's spectral hull lies in the spectrum of its image, giving the reverse of the inclusion in Lemma~\ref{lem:exterior}\ref{item:hull}.

\begin{lemma}\label{lem:periphery}
Let $C$ be a unital $C^*$-algebra, let $D$ be a unital semisimple Banach algebra, and let $S:C\to D$ be a unital linear bijection such that $S(e^xe^y)\in D^\times$ for all $x,y\in C$. For every $x\in C$,
\[
    \partial\widehat{\sigma}_C(x)\subseteq\sigma_D(S(x)) \qquad \text{and}
    \qquad \widehat{\sigma}_D(S(x))=\widehat{\sigma}_C(x).
\]
\end{lemma}

\begin{proof}
By Lemma~\ref{lem:translates}, $S$ is bounded. Suppose that $v\in C$ is the norm limit of a sequence
$(v_n)_{n\in\N}$ of exponentials and that $S(v)\in D^\times$. For each $y\in C$ and $n\in\N$, the left translate gives
\[
    r_D(S(y))
    \stackrel{\eqref{eq:normalized}}{=}
    r_D\left(S(v_n)^{-1}S(v_ny)\right)
    \le \|S(v_n)^{-1}\|\|S\|\|v_ny\|.
\]
By continuity of inversion and multiplication, passage to the limit gives
\begin{equation}\label{eq:approximate-kernel}
    r_D(S(y))\le\|S(v)^{-1}\|\|S\|\|vy\|,
    \qquad y\in C.
\end{equation}

Under the same hypotheses on $v$, suppose that $vh=0$ for some $h\in C$. The inverse $S^{-1}$ is bounded by the open mapping theorem. For $b\in D$, put $y:=S^{-1}(b)$. Since $vh=0$, we have $vy=vS^{-1}(b-S(h))$, so
\begin{align*}
    r_D(b)
    &\stackrel{\eqref{eq:approximate-kernel}}{\le}
      \|S(v)^{-1}\|\|S\|\|vS^{-1}(b-S(h))\|\\
    &\le \|S(v)^{-1}\|\|S\|\|v\|\|S^{-1}\|\|b-S(h)\|.
\end{align*}
By Zem\'anek's characterization~\cite[Theorem~5.3.1(v)]{Aupetit1991}, $S(h)$ belongs to the Jacobson radical of $D$. Since $D$ is semisimple and $S$ is injective, $h=0$.

Now fix $x\in C$ and $\lambda\in\partial\widehat{\sigma}_C(x)$. The boundary of the hull lies in $\sigma_C(x)$, since the hull adds only open complementary components. Thus $v:=\lambda1_C-x$ is singular. Choose a sequence $(\lambda_n)_{n\in\N}$ in $\C\setminus\widehat{\sigma}_C(x)$ converging to $\lambda$. The elements $v_n:=\lambda_n1_C-x\in C$, $n\in\N$, are exponentials by Lemma~\ref{lem:exterior}\ref{item:exterior}. If $v^*v$ were invertible, continuity of inversion would give
\[
    v_n^{-1}=(v_n^*v_n)^{-1}v_n^*
    \longrightarrow (v^*v)^{-1}v^*\qquad\text{in norm}.
\]
Since $v_n\to v$ in norm, passing to the limit in the two inverse identities would give a two-sided inverse of $v$. Hence $0\in\sigma_C(v^*v)$.

For each $n\in \N$, choose a continuous function $f_n:[0,\infty)\to[0,1]$ that equals $1$ on $[0,1/(2n)]$ and vanishes on $[1/n,\infty)$. In $C$ put
\[
    k_n:=f_n(v^*v)\qquad \text{and} \qquad h_n:=f_{2n}(v^*v).
\]
The element $k_n$ is a positive contraction, and $h_n\ne0$ since $0\in\sigma_C(v^*v)$ and $f_{2n}(0)=1$. Moreover, $(1_C-k_n)h_n=0$ because $(1-f_n)f_{2n}=0$ on $[0,\infty)$. Since $1_C-k_n$ is a norm limit of positive invertibles and hence of exponentials, the first part of the proof gives $1_D-S(k_n)\notin D^\times$. Thus $1\in\sigma_D(S(k_n))$. On the other hand, continuous functional calculus gives
\[
    \|vk_n\|^2=\|k_nv^*vk_n\|
    =\max_{t\in\sigma_C(v^*v)}t f_n(t)^2\le\frac1n.
\]
If $S(v)$ were invertible, \eqref{eq:approximate-kernel} with $y=k_n$ would imply $r_D(S(k_n))\to0$, a contradiction.
Thus $\lambda\in\sigma_D(S(x))$.

The unbounded component of $\C\setminus\sigma_D(S(x))$ is disjoint from $\partial\widehat{\sigma}_C(x)$. Since it is
connected and contains all sufficiently large scalars, it lies outside $\widehat{\sigma}_C(x)$. This gives
$\widehat{\sigma}_C(x)\subseteq\widehat{\sigma}_D(S(x))$. The reverse inclusion follows from Lemma~\ref{lem:exterior}\ref{item:hull}.
\end{proof}

\section{Proof of Theorem~\ref{thm:main}}\label{sec:main}

We first prove part~\ref{item:main-single} for two $C^*$-algebras by deriving a resolvent estimate from spectral hull inclusion and using it to show that the second adjoint of the surjection sends projections to idempotents. We then use the results of Section~\ref{sec:mixed} to reduce part~\ref{item:main-two} to part~\ref{item:main-single}.

\subsection{A resolvent estimate}

Let $A$ and $B$ be unital $C^*$-algebras and let $\phi:A\to B$ be a unital linear surjection satisfying the hull inclusion in Theorem~\ref{thm:main}\ref{item:main-single}. We keep these hypotheses through the end of the proof of Proposition~\ref{prop:projections}. By Lemma~\ref{lem:exterior}\ref{item:hull}, $\phi$ is bounded, so the open mapping theorem gives a constant $\kappa>0$ such that every $y\in B$ has a preimage $x\in A$ with $\|x\|\le\kappa\|y\|$.

For a normal element $x_0\in A$, we use this lifting estimate and spectral hull inclusion to bound the resolvent of $\phi(x_0)$ at any $\lambda\in\C$ for which a ray from $\lambda$ stays a positive distance from $\sigma_A(x_0)$.

For $\lambda,\eta\in\C$ with $|\eta|=1$, put
\[
    L(\lambda,\eta):=\{\lambda+\alpha\eta:\alpha\in[0,\infty)\}.
\]
If $x\in A$ and $L(\lambda,\eta)\cap\sigma_A(x)=\emptyset$, the ray is an unbounded connected subset of $\C\setminus\sigma_A(x)$. It follows that $\lambda\notin\widehat{\sigma}_A(x)$, and the hull inclusion gives $\lambda1_B-\phi(x)\in B^\times$.

For nonempty subsets $K_1,K_2\subseteq\C$ and $\lambda\in\C$, we write
\[
    \dist(K_1,K_2):=\inf_{\mu\in K_1,\,\zeta\in K_2}|\mu-\zeta| \qquad \text{and}
    \qquad \dist(\lambda,K_2):=\dist(\{\lambda\},K_2).
\]

\begin{lemma}\label{lem:res}
Let $x_0\in A$ be normal and let $\lambda,\eta\in\C$ with $|\eta|=1$. If $\delta\in(0,\infty)$ satisfies
\[
    \dist(L(\lambda,\eta),\sigma_A(x_0))\ge\delta,
\]
then $\lambda1-\phi(x_0)\in B^\times$ and
\begin{equation}\label{eq:res}
    \|(\lambda1-\phi(x_0))^{-1}\|\le\frac{\kappa}{\delta}.
\end{equation}
\end{lemma}

\begin{proof}
We use the perturbation argument of Colbrook and Stepaniants~\cite[proof of Lemma~3.2]{ColbrookStepaniants2026}. Let $y\in B$ satisfy $\|y\|<\delta/\kappa$ and choose $x\in A$ with $\phi(x)=y$ and $\|x\|\le\kappa\|y\|<\delta$. Fix
$\mu\in L(\lambda,\eta)$. Then $\mu1-x_0$ is invertible. Since $x_0$ is normal, continuous functional calculus gives
\[
    \|(\mu1-x_0)^{-1}\|
    =\dist(\mu,\sigma_A(x_0))^{-1}\le\delta^{-1}.
\]
Consequently $\|(\mu1-x_0)^{-1}x\|<1$. Both factors in
\[
    \mu1-(x_0+x)=(\mu1-x_0)\left(1-(\mu1-x_0)^{-1}x\right)
\]
are invertible, the second by the Neumann series. Since $\mu\in L(\lambda,\eta)$ was arbitrary, the ray is disjoint from $\sigma_A(x_0+x)$. The preceding observation gives $\lambda1-\phi(x_0)-y\in B^\times$.

Thus the open ball of radius $\delta/\kappa$ about $c:=\lambda1-\phi(x_0)\in B$ lies in $B^\times$. In particular,
$c$ is invertible and its distance from $B\setminus B^\times$ is at least $\delta/\kappa$. In a $C^*$-algebra this distance equals $\|c^{-1}\|^{-1}$ by~\cite[Lemma~2.2]{ColbrookStepaniants2026}. This proves~\eqref{eq:res}.
\end{proof}

\subsection{Images of projections in the bidual}

For each projection $p\in A^{**}$, we combine the resolvent estimate in Lemma~\ref{lem:res} with holomorphic functional calculus to show that the second adjoint $\phi^{**}:A^{**}\to B^{**}$ sends $1-2p$ to an element whose square is $1$.

\begin{proposition}\label{prop:projections}
The second adjoint $\phi^{**}:A^{**}\to B^{**}$ sends projections to idempotents.
\end{proposition}

\begin{proof}
Let $p\in A^{**}$ be a projection and put $w:=\phi^{**}(1-2p)\in B^{**}$.
Since $\phi^{**}$ is unital, $\phi^{**}(p)=(1-w)/2$, so it suffices
to show that $w^2=1$.

By Kaplansky's density theorem, there is a net $(k_j)_{j \in \Lambda}$ of self-adjoint contractions in $A$ converging strong$^*$ to $2p-1$. Thus the positive contractions $h_j:=(1+k_j)/2\in A$, $j \in \Lambda$, converge strong$^*$ to $p$. Put
\[
    E:=\{-1,1\} \qquad \text{and} \qquad
    \Gamma_\varepsilon:=\{\zeta\in\C:|\zeta|=1,\
    \varepsilon\operatorname{Im}\zeta\ge 0\},
    \qquad \varepsilon\in E.
\]
For $j \in \Lambda$ and $\varepsilon\in E$, set $u_{j,\varepsilon}:=e^{i\varepsilon\pi h_j}\in A$. Since $\sigma_A(h_j)\subseteq[0,1]$,  the continuous functional calculus shows that $u_{j,\varepsilon}$ is unitary and has spectrum in $\Gamma_\varepsilon$. For each $\varepsilon\in E$, applying Lemma~\ref{lem:functional} with $I:=[0,1]$ and $f(t):=e^{i\varepsilon\pi t}$ for $t\in I$ gives
\[
    u_{j,\varepsilon}\longrightarrow e^{i\varepsilon\pi p}
    =(1-p)+e^{i\varepsilon\pi}p=1-2p
    \qquad\text{in the strong$^*$ topology}.
\]

Fix $\lambda\in\C\setminus E$ and choose $\varepsilon\in E$ so that $\varepsilon\operatorname{Im}\lambda\le 0$. Let $\mu\in L(\lambda,-i\varepsilon)$ and $\zeta\in\Gamma_\varepsilon$. Then $\operatorname{Re}\mu=\operatorname{Re}\lambda$ and $|\operatorname{Im}\mu|\ge|\operatorname{Im}\lambda|$, so $|\mu|\ge|\lambda|$. Since $\operatorname{Im}\mu\,\operatorname{Im}\zeta\le 0$ and $|\operatorname{Re}\zeta|\le1$, we obtain
\[
    |\mu-\zeta|^2=|\mu|^2+1-2\operatorname{Re}(\mu\overline\zeta)
    \ge|\lambda|^2+1-2|\operatorname{Re}\lambda|
    =\dist(\lambda,E)^2.
\]
Set $\delta:=\dist(\lambda,E)>0$. The ray $L(\lambda,-i\varepsilon)$ is at distance at least $\delta$ from $\Gamma_\varepsilon$ and hence from $\sigma_A(u_{j,\varepsilon})$ for each $j \in \Lambda$. Thus Lemma~\ref{lem:res} shows that $c_j:=\lambda1-\phi(u_{j,\varepsilon})\in B$ is invertible for each $j \in \Lambda$. We pass the bound to the bidual as
in~\cite[proof of Lemma~3.2]{ColbrookStepaniants2026}. The $C^*$-identity and continuous functional calculus give
\[
    (\delta/\kappa)^2 1
    \stackrel{\eqref{eq:res}}{\le} c_j^*c_j \qquad \text{and} \qquad
    (\delta/\kappa)^2 1
    \stackrel{\eqref{eq:res}}{\le} c_jc_j^*,
    \qquad j \in \Lambda.
\]
Indeed, the inverses of $c_j^*c_j$ and $c_jc_j^*$ both have norm $\|c_j^{-1}\|^2$.

The net $(c_j)_{j \in \Lambda}$ is bounded since the $u_{j,\varepsilon}$ are unitary and $\phi$ is bounded. By strong$^*$ continuity of $\phi^{**}$ on bounded sets, the net converges strong$^*$ to $c:=\lambda1-w\in B^{**}$. The two products therefore converge strongly to $c^*c$ and $cc^*$, respectively. Testing against vectors in a faithful normal representation of
$B^{**}$ gives the lower bound $(\delta/\kappa)^2 1$ for both limit products, so they are invertible by continuous functional calculus. Thus $c$ has a left inverse $(c^*c)^{-1}c^*$ and a right inverse $c^*(cc^*)^{-1}$, and is therefore invertible. The $C^*$-identity and the lower bound for $c^*c$ give
\begin{equation}\label{eq:limit}
    \|(\lambda1-w)^{-1}\|=\|(c^*c)^{-1}\|^{1/2}
    \le\frac{\kappa}{\dist(\lambda,E)},
    \qquad \lambda\in\C\setminus E.
\end{equation}
As in the proof of~\cite[Lemma~3.3]{ColbrookStepaniants2026}, we use holomorphic functional calculus. Since $\sigma_{B^{**}}(w)\subseteq E$, for $\rho\in(0,1)$ we have
\[
    w^2-1=\frac{1}{2\pi i}
    \sum_{\varepsilon\in E}\oint_{|\lambda-\varepsilon|=\rho}
         (\lambda^2-1)(\lambda1-w)^{-1}\,d\lambda,
\]
where the circles are positively oriented. For $\varepsilon\in E$ and $\lambda\in\C$ with $|\lambda-\varepsilon|=\rho$, we have $\dist(\lambda,E)=\rho$, so
\[
    \|(\lambda^2-1)(\lambda1-w)^{-1}\|
    \stackrel{\eqref{eq:limit}}{\le}\kappa|\lambda+\varepsilon|
    \le\kappa(2+\rho).
\]
It follows that $\|w^2-1\|\le 2\kappa\rho(2+\rho)$. Letting $\rho$ tend to zero gives $w^2=1$.
\end{proof}

\begin{proof}[Proof of Theorem~\ref{thm:main}\ref{item:main-single}]
Let $A$ and $B$ be unital $C^*$-algebras, and let $\phi:A\to B$ be a unital linear surjection. Suppose first that $\phi$ satisfies the hull inclusion, so $\phi$ is bounded by Lemma~\ref{lem:exterior}\ref{item:hull}. Proposition~\ref{prop:projections} and Lemma~\ref{lem:jordan}\ref{item:projection} show that $\phi^{**}$ is a Jordan homomorphism. Restricting to $A$ gives the square identity for $\phi$. A unital Jordan homomorphism sends exponentials to invertible elements by Lemma~\ref{lem:jordan}\ref{item:inverse}. Finally, a map sending exponentials to invertible elements satisfies the hull inclusion by Lemma~\ref{lem:exterior}\ref{item:hull}.
\end{proof}

\begin{remark}
Let $A$ be a unital $C^*$-algebra of real rank zero, let $B$ be a unital semisimple Banach algebra, and let $\phi:A\to B$
be a unital linear surjection sending exponentials to invertible elements. Lemma~\ref{lem:exterior}\ref{item:hull} and the theorem of Cui and Hou~\cite[Theorem~3.1]{CuiHou2004Spectral} imply that $\phi$ is a Jordan homomorphism.
\end{remark}

\subsection{Proof of the two-exponential assertion}

We apply part~\ref{item:main-single} to the bijections obtained from the normalized translates, either on the codomain or, after conjugation, on the quotient. Differentiating the resulting square identities and applying Sinclair's theorem~\cite[Theorem~3.3]{Sinclair1970} gives the derivations needed in Lemma~\ref{lem:tangent}.

\begin{proof}[Proof of Theorem~\ref{thm:main}\ref{item:main-two}]
Let $A$, $B$ and $\phi:A\to B$ be as in Theorem~\ref{thm:main}\ref{item:main-two}, and suppose that $\phi$ satisfies the two-exponential condition. By Lemma~\ref{lem:translates}, $K:=\ker\phi$ is a closed two-sided ideal of $A$. It is proper because $\phi$ is unital. Put $C:=A/K$ and let $\pi:A\to C$ be the quotient map. Define $S:C\to B$ by $S(x+K):=\phi(x)$ for $x\in A$. Then $S$ is a unital linear bijection, and the quotient norm gives $\|S\|\le\|\phi\|$. Its inverse is bounded by the open mapping theorem.

Since $\pi$ commutes with exponentiation, every product of two exponentials in $C$ lifts to such a product in $A$. Hence $S$ satisfies the hypothesis of Lemma~\ref{lem:translates}. If $A$ is a $C^*$-algebra, then its closed two-sided ideal $K$ is self-adjoint, and $C$ is a $C^*$-algebra with its quotient norm.

Fix $a\in C$. For $x\in C$ and $b\in B$, let $L_x:C\to C$ and $L_b:B\to B$ denote left multiplication. The two-exponential condition gives $S(e^{\zeta a})\in B^\times$ for every $\zeta\in\C$. Define $U_\zeta:B\to B$ by
\[
    U_\zeta:=L_{S(e^{\zeta a})^{-1}}SL_{e^{\zeta a}}S^{-1},
    \qquad \zeta\in\C.
\]
Both multiplication operators are invertible and $U_\zeta(1_B)=1_B$, so $U_\zeta$ is a unital linear bijection. Since $U_\zeta S$ is the normalized left translate of $S$ by $e^{\zeta a}$, \eqref{eq:normalized} shows that $U_\zeta$ preserves spectral hulls. Analyticity of inversion and the bounded linear dependence of left multiplication on its multiplier show that the family is entire in operator norm. Since $S$ is unital, we have $U_0=\operatorname{id}_B$, and the derivative of $\zeta\mapsto S(e^{\zeta a})^{-1}$ at zero is $-S(a)$. The product rule gives
\begin{equation}\label{eq:translate-derivative}
    U'_0=SL_aS^{-1}-L_{S(a)}.
\end{equation}

Suppose that $B$ is a $C^*$-algebra. By part~\ref{item:main-single}, every $U_\zeta$ is a Jordan automorphism. Differentiating $U_\zeta(y^2)=U_\zeta(y)^2$ at $\zeta=0$ gives
\[
    U'_0(y^2)=U'_0(y)y+yU'_0(y),\qquad y\in B.
\]
Thus $U'_0$ is a bounded Jordan derivation. Since $B$ is semisimple, it is a derivation by Sinclair's theorem~\cite[Theorem~3.3]{Sinclair1970}. Since $a\in C$ was arbitrary, \eqref{eq:translate-derivative} and Lemma~\ref{lem:tangent} show that $S$ is a Jordan homomorphism.

Suppose now that $A$ is a $C^*$-algebra. By Lemma~\ref{lem:periphery}, $S$ and $S^{-1}$ preserve spectral hulls. For $\zeta\in\C$, the map $V_\zeta:=S^{-1}U_\zeta S:C\to C$ is therefore a unital linear bijection preserving spectral hulls. Part~\ref{item:main-single} shows that every $V_\zeta$ is a Jordan automorphism. Since $V_0=\operatorname{id}_C$, differentiating the square identity
at $\zeta=0$ shows that $V'_0$ is a bounded Jordan derivation on the semisimple algebra $C$, hence a derivation by Sinclair's theorem. From the definition of $V_\zeta$ we obtain
\[
    -V'_0=-S^{-1}U'_0S
    \stackrel{\eqref{eq:translate-derivative}}{=}
    S^{-1}L_{S(a)}S-L_a.
\]
As $a$ ranges over $C$, $S(a)$ ranges over $B$, so the derivations $-V'_0$ satisfy the hypothesis of Lemma~\ref{lem:tangent} for $S^{-1}$. Thus $S^{-1}$, and hence $S$, is a Jordan isomorphism.

In either case, $S\pi=\phi$ is a Jordan homomorphism. Conversely, a unital Jordan homomorphism preserves invertibility by Lemma~\ref{lem:jordan}\ref{item:inverse}, and products of exponentials are invertible.
\end{proof}

\begin{remark}\label{rem:harris-kadison}
Let $A$ and $B$ be as in Theorem~\ref{thm:main}\ref{item:main-two}, and let $\phi:A\to B$ be a linear surjection satisfying the two-exponential condition. Then $c:=\phi(1_A)\in B^\times$, and applying the theorem to $J:=c^{-1}\phi$ gives $\phi=cJ$, where $J:A\to B$ is a unital Jordan homomorphism. When both algebras are $C^*$-algebras, the same conclusion holds if $\phi(e^x)\in B^\times$ for every $x\in A$.

If $A$ and $B$ are $C^*$-algebras and $\phi$ is bijective, then $J$ is a Jordan isomorphism. Thus the individual-exponential condition suffices for the factorization conjectured by Harris and Kadison~\cite[Conjecture~3.1]{HarrisKadison1996}. Under their hypothesis $\phi(A^\times)=B^\times$, the unital bijection $J$ preserves spectra, since $J(A^\times)=c^{-1}B^\times=B^\times$. The conjecture therefore also follows from the spectral isometry theorem of Colbrook and Stepaniants~\cite[Theorem~1.1]{ColbrookStepaniants2026}.
\end{remark}

\subsection*{Funding}
The first author was supported by the European Union -- NextGenerationEU through the National Recovery and Resilience Plan 2021--2026, Institutional Grants of the University of Zagreb Faculty of Science (IK IA 1.1.3. Impact4Math, PMF-CROFUND).

The second author was supported in part by the Slovenian Research Agency (research program P1-0285 and
research projects J1-50000, N1-0296, N1-0428, J1-70047, and J1-70046).

\subsection*{AI use disclosure}
Anthropic's Claude Opus~5.5 found the original proof of the individual-exponential assertion for two $C^*$-algebras. The proof adapts a resolvent estimate of Colbrook and Stepaniants~\cite{ColbrookStepaniants2026}. OpenAI's GPT-6.1 obtained the mixed-algebra extension using normalized translates and assisted with the literature search, proofreading and copyediting. Claude Opus~5.5 also assisted with a subsequent revision and checks of references and selected identities. The initial mixed-algebra proof used the spectral-isometry theorem of Colbrook and Stepaniants, while the present proof uses the spectral hull criterion established in Section~\ref{sec:main}.  The authors take full responsibility for the final content of the manuscript.

\end{document}